\documentclass[12pt]{article}

\usepackage{latexsym,amssymb,amsmath,amsthm,euscript,graphicx}
\usepackage{arydshln}
\usepackage{hyperref}
\usepackage{enumitem}

\newcommand{\DEF}[1]{\textit{#1}}

\newcommand{\aA}{\EuScript{A}}
\newcommand{\bB}{\EuScript{B}}
\newcommand{\DD}{\EuScript{D}}

\newtheorem{thm}{Theorem}
\newtheorem{theorem}[thm]{Theorem}

\newtheorem{lem}[thm]{Lemma}

\newcommand\numberthis{\addtocounter{equation}{1}\tag{\theequation}}

\def\barroman#1{\sbox0{#1}\dimen0=\dimexpr\wd0+1pt\relax
  \makebox[\dimen0]{\rlap{\vrule width\dimen0 height 0.06ex depth 0.06ex}%
    \rlap{\vrule width\dimen0 height\dimexpr\ht0+0.03ex\relax 
            depth\dimexpr-\ht0+0.09ex\relax}%
    \kern.5pt#1\kern.5pt}}

\begin{document}

%~1 1/4 spacing
\setlength{\baselineskip}{1.25 \baselineskip}

%% Title and other pre-body stuff
%%
%% DATE in form: 1 June 1963
%%
\def \today {\number \day \ \ifcase \month \or January\or February\or
  March\or April\or May\or June\or July\or August\or
  September\or October\or November\or December\fi\
  \number \year}

%% TITLE..AUTHOR..DATE
%%

\title{\sffamily Uniquely {\boldmath$D$}-colourable digraphs with large girth \barroman{II}: simplification via generalization} 

\author{
   \textsc{P. Mark Kayll\footnotemark} \\[0.25em]
    {\small\textit{Department of Mathematical Sciences}} \\[-0.25em]
    {\small\textit{University of Montana}} \\[-0.25em]
    {\small\textit{Missoula MT 59812, USA}} \\[-0.1em]
    {\small\texttt{mark.kayll@umontana.edu}}
    \and\addtocounter{footnote}{-1}
    \textsc{Esmaeil Parsa\footnotemark~\footnotemark} \\[0.25em]
    {\small\textit{Department of Mathematics}} \\[-0.25em]
   {\small\textit{Islamic Azad University (Parand Branch)}} \\[-0.25em]
   {\small\textit{Parand New Town, Iran}} \\[-0.1em]
   {\small\texttt{esmaeil.parsa@yahoo.com}}
}

%% \date{}
\date{\small \today}

\maketitle

%
% switch to manual footnotes to avoid pre-superscript on AMS
% Subject. . .
\renewcommand{\thefootnote}{}

%
% Use symbols instead of numbers for footnotes
%
\renewcommand{\thefootnote}{\fnsymbol{footnote}}

\addtocounter{footnote}{1}
\footnotetext{Partially supported by a grant from the Simons Foundation (\#279367 to Mark Kayll)} 
\addtocounter{footnote}{1}
\footnotetext{Partially supported by a 2017 University of Montana
  Graduate Student Summer Research Award funded by the George and Dorothy Bryan
Endowment}
\footnotetext{This work forms part of the author's PhD dissertation~\cite{Parsa}.}

  \begin{abstract}

\noindent
We prove that for every digraph $D$ and every choice of positive
integers $k$, $\ell$ there exists a digraph $D^*$ with girth at least
$\ell$ together with a surjective acyclic homomorphism 
$\psi\colon D^*\to D$ such that: (i) for every digraph $C$ of order at
most $k$,
there exists an acyclic homomorphism $D^*\to C$ if and only
if there exists an acyclic homomorphism $D\to C$; and (ii)
for every $D$-pointed digraph $C$ of order at most $k$ and
every acyclic homomorphism $\varphi\colon D^*\to C$ there exists a unique
acyclic homomorphism $f\colon D\to C$ such that $\varphi=f\circ\psi$. 
This implies the main results in [A. Harutyunyan et al., Uniquely
$D$-colourable digraphs with large girth, \textit{Canad.\ J. Math.},
\textbf{64(6)} (2012), 1310--1328; MR2994666] analogously with how the work 
[J. Ne\v{s}et\v{r}il and X. Zhu, 
On sparse graphs with given colorings and homomorphisms,
\textit{J. Combin. Theory Ser. B}, \textbf{90(1)} (2004), 161--172; MR2041324]
generalizes and extends
[X. Zhu, Uniquely {$H$}-colorable graphs with large girth, 
\textit{J. Graph Theory}, \textbf{23(1)} (1996), 33--41; MR1402136].


   \end{abstract}

\renewcommand{\thefootnote}{\arabic{footnote}}
\addtocounter{footnote}{-2}

\section{Introduction}
\label{S:1}

In 1959, Paul Erd\H{o}s, in a landmark paper~\cite{Erdos1959}---now
known as one of the most pleasing uses of the probabilistic 
method---proved the existence of graphs with arbitrarily large girth and
chromatic number. His technique has been extended in a number of
ways, e.g., 
by Bollob\'as and Sauer~\cite{Bollobas1976} to prove that for all
$k\geq 2$ and $\ell\geq 3$ there is a uniquely $k$-colourable graph whose
girth is at least $\ell$.  It would be difficult to overstate the
influence of this one~\cite{Erdos1959} of Erd\H{o}s' thousands of
results. Indeed, one authoritative combinatorialist
went so far as to assert that 
``All interesting combinatorics flows from the existence of graphs
with large girth and chromatic 
number.''\footnote{St\'{e}phan Thomass\'{e} included the assertion in his
  plenary CanaDAM lecture, 2 June 2011, Victoria, Canada.}
Of course, we interpret Thomass\'{e}'s remark as somewhat
tongue-in-cheek, but as they say, \textsl{many a truth is said in jest.}
In the present article, we follow the flow, from colourings to
homomorphisms and from graphs to digraphs. This work is a sequel
to \cite{Ararat12}, with which we assume some familiarity. For
example, because the introduction of \cite{Ararat12} is more extensive than
this one, we refer the reader there for more background. Also, 
some of the arguments from \cite{Ararat12}---e.g.\ the statement/proof 
of Lemma~\ref{clm3} and Lemma~\ref{boundonPr} (both 
below)---prove useful here. We try to balance
the conflicting goals of not duplicating earlier work while allowing
our new results to stand on their own.

Erd\H{o}s' argument in \cite{Erdos1959} was probabilistic, hence
nonconstructive. To help answer the question of what graphs with large
girth and chromatic number actually look like, in 1968
Lov\'{a}sz~\cite{Lovasz} constructed hypergraphs with arbitrarily large
girth and chromatic number. M\"{u}ller~\cite{muller} also worked in
this domain. More than twenty years after Lov\'{a}sz's contribution,
K\v{r}\'{i}\v{z}~\cite{kriz} produced the first purely graph-theoretic
construction of graphs with arbitrarily large girth and chromatic
number. And more recently (2016), Alon et al.\,\cite{AKRWZ2016} 
constructed such graphs that also satisfy a side condition on
maximum average degree. The time intervals separating these results
offer some hint of the delicacy of their constructions.

Graph homomorphisms, as vertex mappings that preserve
adjacency, naturally generalize graph colouring.
In 1996, working in this realm, Zhu~\cite{X.Zhu1996} 
 proved that for every `core' graph $H$ and every positive integer $\ell\geq 3$ there
exists a uniquely $H$-colourable graph with girth at least $\ell$. 
Because complete graphs are cores, Zhu's result
generalized \cite{Bollobas1976} and \cite{Erdos1959}. 
Almost ten years later, Ne\v{s}et\v{r}il 
and Zhu~\cite{Nesetril2004} further generalized the 
results in the sequence \cite{Erdos1959,Bollobas1976, X.Zhu1996}
using the notion of `pointed' graphs. 

Let us shift now to digraphs. Their circular chromatic number
was first studied in \cite{Mark04}, where Bokal et al.\ 
showed that the colouring theory for digraphs is similar
to that for undirected graphs when stable vertex sets are replaced by
acyclic sets. For example, using an analogue of Erd\H{o}s' original
argument from \cite{Erdos1959} , they showed that there exist 
digraphs of arbitrarily large (directed) girth and circular chromatic
number. Almost a decade later, in \cite{Ararat12}, a subset of these
authors together with their doctoral students
established analogues of Zhu's results from
\cite{X.Zhu1996} in a digraph setting; namely, for a suitable digraph
$D$, there exist digraphs of arbitrarily large girth that are uniquely
$D$-colourable. Severino~\cite{severino1}  presented a construction of
highly chromatic digraphs without short cycles and another
construction \cite{severino2} of uniquely $n$-colourable digraphs (for
arbitrary $n$) with
arbitrarily large girth. The latter two articles, based 
on \cite{Severino-thesis}, give constructive
proofs of results in \cite{Mark04} and \cite{Ararat12} that were
originally proved probabilistically.

This paper analogizes the results of Ne\v{s}et\v{r}il 
and Zhu~\cite{Nesetril2004} to the realm of digraphs.
Just as \cite{Nesetril2004} puts the final icing on the sequence 
\cite{Erdos1959,Bollobas1976,X.Zhu1996}, so too does our main 
result---Theorem~\ref{main} below---provide a fitting capstone for the
sequence \cite{Mark04,Ararat12}. Postponing definitions for another
minute (until Section~\ref{Background-terminology}), let us state our
main result and lay bare its connection with \cite{Ararat12}.

\begin{theorem}
\label{main}
 For every digraph $D$ and every choice of positive integers $k$, 
 $\ell$ there exists a digraph $D^*$ together with a surjective acyclic
 homomorphism $\psi\colon D^*\to D$ with the following properties: 
\begin{enumerate}
\renewcommand{\theenumi}{\upshape(\roman{enumi})}
    \item $\text{girth}(D^*)\geq\ell$;
    \item for every digraph $C$ with at most $k$ vertices, there
      exists an acyclic homomorphism $D^*\to C$ if and only if
      there exists an acyclic homomorphism $D\to C$;
    \item for every $D$-pointed digraph $C$ with at most $k$ vertices
      and for every acyclic homomorphism $\varphi\colon D^*\to C$ there
      exists a unique acyclic homomorphism $f\colon D\to C$ such
      that $\varphi=f\circ\psi$.  
\end{enumerate}
\end{theorem}

The precursor~\cite{Ararat12} established two main results:

\begin{theorem}
\label{first}
 If $D$ and $C$ are digraphs such that $D$ is not $C$-colourable, then
for every positive integer $\ell$, there exists a digraph $D^*$ of girth at
least $\ell$ that is $D$-colourable but not $C$-colourable. 
\end{theorem}

\begin{theorem}
\label{second}
 For every core $D$ and every positive integer $\ell$, there is a
 digraph $D^*$ of girth at least $\ell$ that is uniquely $D$-colourable.
\end{theorem}

To see that Theorem~\ref{main} implies Theorem~\ref{first}, let us be
given a positive integer $\ell$ and two digraphs $C$, $D$ with 
$D$ not $C$-colourable (as in the hypotheses of Theorem~\ref{first}).
Taking $k$ to be the order of $C$, we can put this $C$ in the role of 
the digraph $C$ in conclusion~(ii) of Theorem~\ref{main}, which
delivers a digraph $D^*$ with $D^*\to D$. As
$D\not\rightarrow C$, the same conclusion shows that also 
$D^*\not\rightarrow C$, and conclusion~(i) gives the girth requirement
on $D^*$.

Before deriving Theorem~\ref{second} from Theorem~\ref{main},
observe that if $D$ is a core, then every acyclic homomorphism from
$D$ to itself must be an automorphism, and so if any two such
homomorphisms agree on all but one vertex, they must 
also agree on that vertex. Therefore, cores $D$ are $D$-pointed. 

Now let us be given a positive integer $\ell$ and a core $D$ (as in
the hypotheses of Theorem~\ref{second}). If we here take $k =|V (D)|$,
then Theorem~\ref{main} delivers a large-girth digraph $D^*$ together with a
$D$-colouring $\psi\colon D^*\to D$.  The preceding paragraph foreshadows
that we can put $D$ in the role of $C$ in conclusion~(iii), which
shows that every acyclic homomorphism $\varphi\colon D^*\to D$ 
yields an acyclic homomorphism $f\colon D\to D$ such that
$\varphi=f\circ\psi$.  But $D$ being a core implies that such an $f$
is an automorphism, so we've shown that $\varphi$ and $\psi$ differ by
an automorphism, i.e., that $D^*$ is uniquely $D$-colourable.

Notice that being $D$-pointed is a necessary condition in part~(iii)
of Theorem~\ref{main}. For consider two acyclic homomorphisms
$f',f''\colon D\to C$ satisfying (for some vertex $x_0$ of $D$) 
$f'(x)=f''(x)$ for all $x\neq x_0$
and $f'(x_0)\neq f''(x_0)$, and assume that there is an arc between
$f'(x_0)$ and $f''(x_0)$ in $C$. Typically, the set $\psi^{-1}(x_0)$ can be split
into two nonempty sets $A$, $B$ and we can define $\varphi\colon D^*\to C$
by $f'\circ\psi(y)$ for $y\in V(D^*)\smallsetminus B$ and 
$f''\circ\psi(y)$ for $y\in B$. Now this $\varphi$ sends $A$ and $B$ to two
different points while $f\circ\psi$, for any given $f\colon D\to
C$, sends these sets to a single point. Therefore, the acyclic
homomorphism $\varphi$ cannot be written as $\varphi=f\circ\psi$ 
for an acyclic homomorphism $f\colon D\to C$. 

\subsection*{Remarks}
As hinted above, Ne\v{s}et\v{r}il's and Zhu's
article~\cite{Nesetril2004} was in a sense a crowning achievement for
a body of work initiated by Erd\H{o}s in \cite{Erdos1959}. For any
given graph $G$, they produced a high-girth graph $G^*$ characterizing
the small-order graphs admitting a homomorphism from $G$ and
furthermore, via $G$-pointedness, wound unique colourability into their
tapestry. Their results generalized \cite{Bollobas1976}, \cite{X.Zhu1996} and
moreover some other major contemporary theorems (e.g., the Sparse
Incomparability Lemma and M\"{u}ller's Theorem---see \cite{X.Zhu1996}
and the discussion in \cite{Nesetril2004}).

Because our Theorem~\ref{main} likewise characterizes when
the high directed girth, high digraph chromatic number (for unique
colourability) phenomenon occurs---phrased in terms of acyclic
homo\-morphisms---it too reaches a satisfying destination, now
for the sequence \cite{Mark04,Ararat12}. And because this level of
generality has actually shortened the proofs from \cite{Ararat12},
perhaps we've arrived at the `right' vantage point for viewing these
results.

\section{Terminology, notation, and an auxiliary result}
\label{Background-terminology}

Without being overly encyclopedic, we attempt to include the required
definitions. For basic notation and terminology concerning graphs and
digraphs, we mainly follow \cite{graph} and \cite{digraph},
respectively, and we refer the reader there for any omissions. For a
more (most) thorough treatment of graph homomorphisms, the reader could
consult \cite{Godsil01} (\cite{Hell04}). For probabilistic concerns,
see, e.g., \cite{AlonSpencer_4e-2016} or \cite{MolloyReed-2002}.

All our digraphs are finite and \DEF{simple}---i.e.\ loopless and without 
multiple arcs---however, we do allow two vertices $u$, $v$ to be
joined by two oppositely 
directed arcs $uv$, $vu$. \DEF{Cycles} in digraphs mean directed 
ones, and the \DEF{girth} of a digraph $D$ is the length of a shortest
cycle in $D$.

Just as graph homomorphisms generalize graph colouring, so too do 
acyclic homomorphisms of digraphs generalize (one variant of) digraph 
colouring. So we begin by recalling the definition of these sorts of
homomorphisms from \cite{Mark04}; see \cite{Ararat12} for background.
An \DEF{acyclic homomorphism} of a digraph $D$ to a digraph $C$ is a
function $\rho\colon V(D)\to V(C)$ such that: 
\begin{enumerate}
\renewcommand{\theenumi}{(\roman{enumi})}
    \item for every arc $uv\in A(D)$, either $\rho(u)=\rho(v)$, or
      $\rho(u)\rho(v)$ is an arc of $C$; and  
        \item for every vertex $x\in V(C)$, the subdigraph of $D$
          induced by $\rho^{-1}(x)$ is acyclic. 
\end{enumerate}
Acyclic homomorphisms can also be viewed as a generalization of (ordinary)
homomorphisms of undirected graphs; again, see \cite{Ararat12}.

If there exists an acyclic homomorphism of $D$ to $C$, we say that $D$ is
\DEF{homomorphic} to $C$ and write $D\to C$. 
Motivated by the connection to `acyclic digraph colouring', we
sometimes call an acyclic homomorphism of $D$ to $C$ 
a \DEF{$C$-colouring} of $D$ and say that $D$ is \DEF{$C$-colourable}.
A digraph $D$ is \DEF{uniquely $C$-colourable} if it is surjectively 
$C$-colourable, and for any two $C$-colourings $\psi$, $\varphi$
of $D$, there is an automorphism $f$ of $C$ such that $\varphi=f\circ\psi$;
when this occurs, we say that $\varphi$ and $\psi$ 
\DEF{differ by an automorphism} of $C$. 
A digraph $D$ is a \DEF{core} if
the only acyclic homomorphisms of $D$ to itself are
automorphisms. Given two digraphs $C$, $D$, we say that $C$ is
\DEF{$D$-pointed} if there do 
not exist two $C$-colourings $\rho$, $\varphi$ of $D$
such that $\rho(v)\neq\varphi(v)$ holds for exactly one vertex $v$  of
$D$. As noted following the statement of Theorem~\ref{second}, digraph
cores $D$ are $D$-pointed.

\subsection*{Probabilistic tools}

Our proof of Theorem~\ref{main} invokes several standard probabilistic
tools. Aside from the First Moment Method (Markov's
Inequality)---which is explicitly invoked a handful of
times---Inclusion-Exclusion and the Janson Inequalities also make an
implicit appearance through their use (in \cite{Ararat12}) in proving 
Lemma~\ref{boundonPr} below. We shall not restate these standard 
results here; however, for convenience, we do include a version of 
Chernoff's famous bound(s) on the tail distributions of binomial
random variables. Though more technical versions are 
available---see, e.g., \cite{JansonLuczakRucinski2000}---this
one will suffice for our main proof in Section~\ref{Sec:main_proof}:

\begin{thm}
\label{Chernoffresult}
If $X$ is a binomial random variable and $0<\gamma<3/2$, then
\[
P(|X-E(X)|\geq\gamma E(X))\leq 2 e^{-\gamma^2E(X)/3}.
\]
\end{thm}
\noindent

\section{\protect{Set-up for the proof of Theorem~\ref{main}}}
\label{s:3}

We begin at the starting point for the main proof in \cite{Ararat12},
namely specifying a random digraph model, which needs no change here.
Suppose that the digraph $D$ is given with $V(D)=\{1,2,\ldots,a\}$ and
$|A(D)|=q$. Let $n$ be a positive integer and 
$V_1,V_2,\ldots,V_a$ be pairwise-disjoint ordered $n$-sets
$V_i=\{v_{i_1},v_{i_2},\ldots,v_{i_n}\}$, for $i=1,2,\ldots,a$. 
Next let $D_0$ be the digraph with vertex set 
$V:=V_1\cup V_2\cup \cdots \cup V_a$ and  
\begin{align*}
    A(D_0):= & \left\{\text{\rule[-0.5em]{0pt}{2pt}}xy:x\in V_i, y\in V_j\; \text{with}\; ij\in A(D), \text{for some}\; i,j\in\{1,2,\ldots,a\}\right\} \\ &  \text{\rule[-10pt]{6em}{0pt}}\bigcup_{i=1}^a~\left\{\text{\rule[-0.75em]{0pt}{2pt}}v_{i_k}v_{i_t}:k,t\in\{1,2,\ldots,n\}\; \text{and}\; k<t\right\};
\end{align*}
so $D_0$ has $na$ vertices and $a\binom{n}{2}+qn^2$ arcs.

Now fix an $\epsilon$ with $0<\epsilon <1/4\ell$.  Our random digraph
model $\DD(n,p)$ consists of all spanning subdigraphs of $D_0$ in which
the arcs are chosen randomly and independently with probability
$p:=n^{\epsilon-1}$. Through the following three lemmas we prove
essential technical facts about digraphs in $\DD(n,p)$. Throughout the discussion
$n$ is assumed to be sufficiently large to support our assertions. 

Our first aim is to show that most digraphs in $\DD(n,p)$ have few
short cycles which are pairwise vertex-disjoint. 
\begin{lem}
\label{clm1}
\begin{enumerate}
\renewcommand{\theenumi}{\upshape(\roman{enumi})}
\setlength{\itemsep}{0pt}
  \item\label{clm1-i} The expected number of cycles of length less
    than $\ell$ in a digraph $\widehat{D}\in \DD(n,p)$ is bounded from
    above by $n^{\epsilon \ell}n^{-\epsilon/2}$; 
    \item\label{clm1-ii}  the expected number of pairs of cycles of
      length less than $\ell$ in a digraph $\widehat{D}\in \DD(n,p)$
      which intersect in at least one vertex is bounded from above by
      $n^{-1/2}$. 
\end{enumerate}
\end{lem}
\noindent
By Markov's Inequality, Lemma~\ref{clm1} implies that
asymptotically almost all digraphs from $\DD(n,p)$ have at most
$n^{\epsilon \ell}$ cycles of length less than $\ell$, and these
cycles are all vertex-disjoint. The ideas in the proofs of (i) and (ii)
are contained, respectively, in the ``Proof of (2.1)'' and ``Proof of
(3.1)'' in \cite{Ararat12}; we include the proofs here for context,
completeness, and consolidation. 

\begin{proof}
\ref{clm1-i} Let $\widehat{D}\in \DD(n,p)$ and let the random
variables $X_i$, $X$ count, respectively, the
number of cycles of length $i$, for $2\leq i<\ell$, and of length less than $\ell$
in $\widehat{D}$. Then 
\begin{align*}
    E(X_i) \leq \binom{na}{i}(i-1)!p^i=\frac{na(na-1)\cdots (na-i+1)}{i}p^i  < \frac{(na)^i}{i}p^i.
\end{align*}

\noindent
Hence 
\begin{align*}
    E(X)=\sum_{i=2}^{\ell-1} E(X_i)&\leq \sum_{i=2}^{\ell-1} \frac{(na)^i}{i}p^i 
    \leq \sum_{i=2}^{\ell-1} \frac{(n^{\epsilon}a)^i}{i},
        %\;\;\;\;(\;p=n^{\epsilon-1})\\ &
    %< a^{l-1}n^{(l-1)\epsilon}=a^{l-1}n^{-\epsilon}n^{\epsilon l} (*)\\ &
    %<n^{\epsilon l}n^{-\epsilon/2}\;\;\;\;(n \;\text{is large})
\end{align*}

\noindent
recalling that $p=n^{\epsilon-1}$ for the last step. Now, the inequality $\;\sum_{i=2}^{\ell-1}(n^{\epsilon}a)^{i}/i< a^{\ell-1}n^{(\ell-1)\epsilon}\;$ (which can be proved by induction on $\ell$) shows that 
\begin{align*}
    E(X)<a^{\ell-1}n^{(\ell-1)\epsilon}=a^{\ell-1}n^{-\epsilon}n^{\epsilon\ell}<n^{\epsilon\ell}n^{-\epsilon/2},
\end{align*}
for sufficiently large values of $n$.

To prove part~\ref{clm1-ii}, 
we need the
following definition from~\cite{Ararat12} which in turn had its roots
in~\cite{X.Zhu1996}. For integers $\ell_1,\ell_2<\ell$, we call a
digraph an $(\ell_1,\ell_2)$-\textit{double cycle} if it consists of a
directed cycle $C_{\ell_1}$ of length $\ell_1$ and a directed path of
length $\ell_2$ joining two (not necessarily distinct) vertices of
$C_{\ell_1}$. An $(\ell_1,\ell_2)$-double cycle contains
$\ell_1+\ell_2$ arcs and $\ell_1+\ell_2-1$ vertices. 

A moment's reflection shows that if two cycles of length less than
$\ell$ intersect in at least one vertex, then they contain (as a
subdigraph) an $(\ell_1,\ell_2)$-double cycle for some
$\ell_1,\ell_2<\ell$. Hence in a random $\widehat{D}\in \DD(n,p)$ the
expected number of pairs of cycles of length less than $\ell$ that
intersect in at least one vertex is at most the expected number of all
$(\ell_1,\ell_2)$-double cycles for $\ell_1,\ell_2<\ell$.

Let the random variable $Y$ count the number of all
$(\ell_1,\ell_2)$-double cycles for some $\ell_1,\ell_2<\ell$ in a
random  $\widehat{D}\in \DD(n,p)$. For fixed $\ell_1,\ell_2<\ell$, let
$Y(\ell_1,\ell_2)$ be the number  of $(\ell_1,\ell_2)$-double
cycles. Then 
\begin{align*}
      E(Y(\ell_1,\ell_2)) &< 2\binom{an}{\ell_1}(\ell_1-1)!p^{\ell_1}(\ell_1)(\ell_1)\binom{an}{\ell_2-1}(\ell_2-1)!p^{\ell_2}\\ &
       <\ell_1(na)^{\ell_1}(na)^{\ell_2-1}p^{\ell_1+\ell_2} \\ &
       <\ell_1a^{\ell_1+\ell_2}n^{\epsilon(\ell_1+\ell_2)}n^{-1}.
    \end{align*}
    
\noindent
As $\epsilon(\ell_1+\ell_2)<2\ell\epsilon<1/2$ (because
$\ell_1,\ell_2\leq \ell$ and $\epsilon<1/4\ell$), for large enough $n$
we have 
   \begin{align*}
       & E(Y)=\sum_{\substack{2\leq \ell_1<\ell\\ 1\leq \ell_2<\ell}}E(Y(\ell_1,\ell_2))<n^{-1/2}.
       \end{align*}
\end{proof}

To state the second lemma we need the following definition (which
leans on the parameters $D$ and $k$ of Theorem~\ref{main}). This
set-up and the ensuing analysis in Lemma~\ref{clm2} is modelled after
an analogous discussion in~\cite{Nesetril2004}. 
Following these authors, we call a set $\aA\subseteq V$
\DEF{large} if there are distinct $i,j\in [a]$, with
$ij$ an arc of $D$, such that both $|\aA\cap V_i|\geq n/k$ and
$|\aA\cap V_j|\geq n/k$, and the $D$-arc $ij$ in this case is a
\DEF{good} arc for $\aA$. For a large set $\aA$, denote by $|\widehat{D}/\aA|$ the
minimum number of arcs of (a random) $\widehat{D}$ which lie in a set 
$\{xy: x\in\aA\cap V_i,\; y\in\aA\cap  V_j\}$, with $ij$ a good arc for $\aA$.  
 
 \begin{lem}\label{clm2}
 If $\widehat{D}\in \DD(n,p)$ and $\aA$ is large, then $P(|\widehat{D}/\aA|\geq n)=1-o(1)$.
\end{lem}

\noindent
Thus asymptotically most digraphs from $\DD(n,p)$ enjoy the property
of all good arcs (of $D$) for large sets $\aA$ inducing at least $n$ arcs
(of $\widehat{D}\in \DD(n,p)$). 
\begin{proof}
Let $\widehat{D}\in \DD(n,p)$ and $\aA\subseteq V$ be a  large set and set 
$\alpha=P(|\widehat{D}/\aA|\geq n)$.
Essentially following \cite[Proof (of Claim 2)]{Nesetril2004}, we have 
  \begin{align*}
     1-\alpha=P(|\widehat{D}/\aA|<n)& \leq \sum_{\text{$\bB$ large}}P(|\widehat{D}/\bB|<n)\\ &
     \leq 2^{na}\binom{qn^2}{n}(1-p)^{n^2/k^2-n} \\ &
     <e^{cn\ln n-c'n^{1+\epsilon}}=o(1) \numberthis \label{firstbound}
 \end{align*}
\noindent
for some positive constants $c$ and $c'$ that are independent of $n$
(with the estimates in (\ref{firstbound}) 
being borrowed from \cite{Nesetril2004}). Thus we get $\alpha=1-o(1)$.
 \end{proof} 
 
The last lemma of this section addresses a technical 
situation also encountered at the end of Section~3
of~\cite{Ararat12}. We repeat part of the proof here 
for completeness and also to
facilitate fleshing out more of its details. See also 
\cite[Claim~3]{Nesetril2004} for an analogous statement
(for graphs and homomorphisms) and an alternate proof approach (via
enumeration).

\begin{lem}\label{clm3}
  Almost all digraphs from $\DD(n,p)$ do not contain two nonempty sets
  $\aA\subset V_{i_0}$, $\bB\subset V_{j_0}$, for some $i_0,j_0\in[a]$, with
  $i_0j_0\in A(D)$ (resp. $j_0i_0\in A(D)$), $|\aA|=n-(k-1)|\bB|$,
  $|\bB|\leq n/k$, such that the set $\aA\cup\bB$ contains at most 
  $\min\{|\bB|,n^{\epsilon \ell}\}$ arcs from $\aA$ to $\bB$ (resp. from $\bB$ to
  $\aA$) and these arcs form a matching (i.e.\ a set of independent arcs). 
  \end{lem}

\begin{proof}
Let $b\leq n/k$ and 
$s\leq \min \{b,\lceil n^{\epsilon \ell} \rceil\}$. 
We denote by $L(b,s)$ the expected number of pairs $\aA$, $\bB$ 
such that $\aA\subseteq V_i$, $\bB\subseteq V_j$, $ij\in A(D)$, 
$|\aA|=n-(k-1)|\bB|$, $|\bB|=b$ and there are exactly $s$ arcs 
joining a vertex in $\aA$ to a vertex in $\bB$. Then
 \begin{align}
     L(b,s) &\leq  \binom{n}{n-(k-1)b}\binom{n}{b}\binom{(n-(k-1)b)b}{s}p^s(1-p)^{(n-(k-1)b)b-s}\nonumber
     \\ &
     <n^{(k-1)b}n^b(nb)^sn^{s(\epsilon -1)}e^{-bn^{\epsilon}+n^{\epsilon-1}((k-1)b^2+s)}\nonumber
     \\ &
     <n^{kb}b^sn^{\epsilon s}e^{-(bn^{\epsilon})/2}  \label{5.1}
     \\ &
     =b^sn^{\epsilon s}(n^k e^{-n^{\epsilon}/2})^b \nonumber
     \\ &
     <b^sn^{\epsilon s}e^{-(bn^{\epsilon})/3}  \label{5.2} 
     \\ &
     <e^{-n^{\epsilon}/4}. \label{5.3}
     \end{align}
     
     To help the reader through steps (\ref{5.1})--(\ref{5.3}), we fill in
     the following estimates:\\
\noindent
     for (\ref{5.1}):
\[
      -bn^{\epsilon}+n^{\epsilon-1}((k-1)b^2+s)=-bn^{\epsilon}+\frac{(k-1)b^2+s}{n^{1-\epsilon}} < -bn^{\epsilon}+\frac{bn^{\epsilon}}{2}=-\frac{bn^{\epsilon}}{2};
\]
\noindent     
     for (\ref{5.2}): for large enough $n$, we have $n^k<e^{n^{\epsilon/6}}$, so that $n^ke^{-n^{\epsilon/2}}<e^{-n^{\epsilon}/3}$;\\
\noindent
     and lastly for (\ref{5.3}):
\[
      (bn^{\epsilon})^s=e^{s\ln(bn^{\epsilon})}=(e^{\ln(bn^{\epsilon})})^s
     <(e^{(1/12bn^{\epsilon})^{1/s}})^s=e^{1/12bn^{\epsilon}},
\]
 
\noindent
 and this implies that 
\[
      b^sn^{\epsilon s}e^{-(bn^{\epsilon})/3}<e^{-bn^{\epsilon}/4}<e^{-n^{\epsilon/4}}.
\]

 So with $L(b):=\sum_{s\leq \min \{b,\lceil n^{\epsilon \ell} \rceil\}} L(b,s)$, we find that 
\[
     L(b)< \lceil n^{\epsilon \ell}\rceil e^{-n^{\epsilon}/4}< e^{-n^{\epsilon}/5},
\] 
\noindent
 and we finally obtain 
\[
     \sum_{1\leq b\leq n/k}L(b)< (n/k)e^{-n^{\epsilon}/5}<e^{-n^{\epsilon}/6}.
\]
\noindent
An application of Markov's Inequality completes the proof. (Notice
that we are getting a small upper estimate here even without the
matching condition). 
 \end{proof}
 
\section{\protect{Proof of Theorem~\ref{main}}} 
\label{Sec:main_proof}

We continue to be guided by \cite{Nesetril2004}, but the argument here
is complicated by the more technical definition of 
`acyclic homomorphism' in our context compared to `homomorphism'
in the graph setting.

Choose a digraph $D'$ in $\DD(n,p)$ satisfying the properties
 asserted in Lemmas~\ref{clm1}--\ref{clm3}. So $D'$ contains 
at most $n^{\epsilon \ell}$ (directed) cycles of length
less than $\ell$ and these cycles are pairwise vertex-disjoint. 
Consequently  (picking one arc 
from each cycle), there is a matching (an independent arc set) 
$M\subseteq A(D')$ of size at most $n^{\epsilon \ell}$ such that the digraph 
$D'-M=(V(D'),A(D')\smallsetminus M)$ has no cycles of length less than
$\ell$. We prove that this digraph---henceforth denoted
$D^*:=D'-M$---satisfies the conclusions of Theorem~\ref{main}.  

Define $\psi\colon V(D^*)\to V(D)$ by $\psi(x)=i$ if and only if
$x\in V_i$, for $i\in [a]$. It is clear from the definition of
$\DD(n,p)$ that $\psi$ is a surjective acyclic homomorphism. That
$\text{girth}(D^*)\geq \ell$ was arranged in our description of
$D^*$, and this takes care of (i). 
 
To prove part (ii) of Theorem~\ref{main}, fix a digraph $C$ of order
at most $k$ and consider an acyclic homomorphism
$\varphi\colon D^*\to C$. We proceed to define a mapping 
$f\colon V(D)\to V(C)$. By the Pigeonhole
Principle, for each $i\in V(D)$, there is a vertex $x\in V(C)$
such that $|V_i\cap \varphi^{-1}(x)|\geq n/k$.
We let $f(i)=x$ (choosing $x$ arbitrarily if more than one $x$ 
has this property) and now prove that $f$ is an acyclic
homomorphism. To prove that $f$ satisfies the first property of
being an acyclic homomorphism, let $ij$ be an arc of $D$ with
$f(i)=x$ and $f(j)=y$. If $x=y$, then we are done, so suppose that
$x\neq y$. With 
$A_i=V_i\cap \varphi^{-1}(x)$ and $A_j=V_j\cap\varphi^{-1}(y)$,
we have $|A_i|\geq n/k$ and $|A_j|\geq n/k$ from the
definition of $f$. Hence $\aA=A_i\cup A_j$ is a large set and $ij$ is
a good arc for $\aA$, so we can invoke  Lemma~\ref{clm2} to see
that there exists an arc of $D^*$ from $A_i$ to $A_j$ (Note that we
deleted at most $n^{\epsilon\ell}<n^{1/4}$ arcs from $D'$ to get
$D^*$, but $ij$ induces at least $n$ arcs, so we did not delete all
these arcs from $A_i$ to $A_j$). Now, since $\varphi$ is an acyclic
homomorphism, we have $xy\in A(C)$ as required. 

To finish the proof that $f$ is an acyclic homomorphism, we need to
show that $f^{-1}(x)$ induces an acyclic subdigraph in $D$ for every
$x\in V(C)$. We prove this by contradiction. 
Suppose that there is a vertex $v'\in V(C)$ such that the subdigraph
induced by $f^{-1}(v')$ in $D$ contains a cycle $Q$. Write
$Q=i_1i_2\cdots i_t$ and observe that $2\leq t\leq a$. Since
$f(i_s)=v'$, for $s=1,2,\ldots,t$, we have 
$|V_{i_s}\cap\varphi^{-1}(v')|\geq n/k$, for  $s=1,2,\ldots,t$ (from the
definition of $f$). The fact that $n^{\epsilon \ell}\ll n/k$ implies
that each set $V_{i_s}\cap \varphi^{-1}(v')$ contains a subset
$W_{i_s}$ of size $w:=\lceil n/(2k) \rceil$ such that no arc in $M$
has an end vertex in $W_{i_s}$. It follows from 
$W_{i_s}\subseteq V_{i_s}\cap \varphi^{-1}(v')$ that
$\varphi(W_{i_1})=\cdots=\varphi(W_{i_s})=\{v'\}$. Since $\varphi$ is an
acyclic homomorphism, the subdigraph of $D^*$ induced by 
$W_{i_1}\cup W_{i_2}\cup \cdots \cup W_{i_s}$ is acyclic. We show that the event
that $W_{i_1}\cup W_{i_2}\cup \cdots \cup W_{i_s}$ induces an acyclic
subdigraph in $D^*$ is unlikely.  
 
Let us consider a sequence of sets $U_{j_1},U_{j_2},\ldots, U_{j_r}$
such that for $i=1,2,\ldots,r$ we have $U_{j_i}\subseteq V_{j_i}$ and
$|U_{j_i}|=w$, and the vertex sequence $j_1,j_2,\ldots,j_r$ is a cycle
in $D$. We denote by $P_r$ the probability that the subdigraph of
$D^*$ induced by $U_{j_1}\cup U_{j_2}\cup \cdots \cup U_{j_r}$ is
acyclic and call this sequence \DEF{bad} if it induces
an acyclic subdigraph in $D^*$. Now, for the expected number $N$ of
bad sequences in $D^*$, we have 
\begin{equation}
\label{N-bnd}
   N\leq \sum_{r=2}^a\binom{a}{r}(r-1)!\binom{n}{w}^rP_r.
 \end{equation}
We pause to note that (\ref{N-bnd}) is relation~(2.6) from
\cite{Ararat12}, adapted to our present context. 
The following result bounds the probabilities
$P_r$;  for a proof, see \cite{Ararat12} (which actually
contains two proofs).

 \begin{lem}[\protect{\cite[Lemma~2.1]{Ararat12}}] 
 \label{boundonPr}
For every integer $r\in\{2,...,a\}$, we
   have $P_r\leq e^{-n^{1+\epsilon}/(10k^2)}$. 
\end{lem}

\noindent 
As in \cite[relations (2.19)]{Ararat12}, we see that
for large enough $n$, the relation~(\ref{N-bnd})  and Lemma~\ref{boundonPr} give
\[
     N \leq \sum_{r=2}^a \binom{a}{r}(r-1)!\binom{n}{w}^r e^{-n^{(1+\epsilon)}/(10k^2)}
     < \sum_{r=2}^a \frac{e^{-n}}{2a} <\frac{e^{-n}}{2}.
\]

\noindent
So to finish this chain of reasoning as in \cite{Ararat12},
using Markov's Inequality, we find that 
\[
      P(\exists\; \text{a bad sequence})<\frac{e^{-n}}{2}.
\]
  
This achieves the goal stated before Lemma~\ref{boundonPr}
which in turn contradicts our assumption that the subdigraph induced
by $f^{-1}(v')$ in $D$ contains a cycle. Thus 
the forward implication in part~(ii) of
Theorem~\ref{main} is proved.

For the converse in (ii), let $f\colon V(D)\to V(C)$ be
an acyclic homomorphism. We define a mapping 
$\varphi\colon V(D^*)\to V(C)$ as $\varphi(x)=f(i)$, where $x\in
V_i$. Each $V_i$ induces an 
acyclic set in $D^*$. Arcs of $D$ that are mapped to single vertices
$f(i)$ in $C$ do not lead to cycles in preimages
$\varphi^{-1}(f(i))$ because $f$ is itself an acyclic
homomorphism. Furthermore, each arc $xy\in A(D^*)$ with 
$\varphi(x)\not=\varphi(y)$ is mapped to the arc 
$\varphi(x)\varphi (y)\in A(C)$ again because $f$ is
an acyclic homomorphism. Hence $\varphi$ is an acyclic
homomorphism. This completes the proof of part~(ii). 
  
We turn our attention to part~(iii) of Theorem~\ref{main}. Let $C$
be a $D$-pointed digraph of order at most $k$ and $\varphi$ be an
acyclic homomorphism from $D^*$ to $C$. We want to show that there
exists a unique acyclic homomorphism $f\colon V(D)\to V(C)$ such
that $\varphi=f \circ \psi$. Note that for every $i\in V(D)$ there exists
a unique $x_i\in V(C)$ such that $|\varphi^{-1}(x_i)\cap V_i|\geq n/k$. 
Existence follows from the Pigeonhole Principle. If there were
$x'_i\not=x_i$ with the same property 
($|\varphi^{-1}(x'_i)\cap V_i|\geq n/k$), then our definition of $f$ 
here would lead to another
acyclic homomorphism $f'\colon V(D)\to V(C)$ such that
$f(j)=f'(j)$ for all $j\in V(D)\smallsetminus\{i\}$. But then the
$D$-pointedness of $C$ would force
$x_i=f(i)=f'(i)=x'_i$. Now, we define 
$f\colon V(D)\to V(C)$ as $f(i)=x_i$ for $i=1,2,\ldots,a$. Because $f$ is
defined as in part~(ii), we again see that this function is an
acyclic homomorphism. Hence, it remains to show that 
$\varphi=f \circ \psi$.

\paragraph*{Remark}
Until now, parts of our proof have involved carefully piecing together
ideas from \cite{Ararat12} and \cite{Nesetril2004}. The remainder of
the argument follows quite a different path and underscores the extra
complexity inherent in working with acyclic homomorphisms (of
digraphs) compared to ordinary homomorphisms (of graphs).

\subsection*{Proof of \boldmath{$\varphi=f \circ \psi$}}
First, we show that $\varphi$ and $f \circ \psi$ have the same
range. It is clear that Range$(f \circ \psi)\subseteq$
Range$(\varphi)$. To prove the reverse containment, suppose to the
contrary that there is a  vertex $y\in$ Range$(\varphi)$ that is not
in Range$(f \circ \psi)$. Since $y$ is in the range of $\varphi$, the 
set $\varphi^{-1}(y)\cap V_i$ is not empty for some $i\in \{1,2,\ldots,a\}$. %Furthermore, Since $y$ is not in the range of $f$, $|\varphi^{-1}(y)\cap V_i|<n/k$.
On the other hand, the definition of $f$ shows that
$|\varphi^{-1}(f(i))\cap V_i|\geq n/k$; in particular
$\varphi^{-1}(f(i))\cap V_i \neq \varnothing$. Because 
$f(i)\in\text{Range}(f \circ \psi)$ while 
$y\not\in \text{Range}(f\circ \psi)$ we see that 
$V_i\smallsetminus \big(\varphi^{-1}(f(i))\cap V_i\big)\neq \varnothing$ 
for some $i\in \{1,2,\ldots,a\}$. We show that this leads to a contradiction.

Let $i_0\in \{1,2,\ldots,a\}$ be such that
$t:=|\varphi^{-1}(f(i_0))\cap V_{i_0}|$ is minimum. It is easy to see
that $t\geq n/k$. Our discussion in the preceding paragraph implies that
$t<n$. We choose $x\in V(C)$ with $x\neq f(i_0)$ such that
$b:=|\varphi^{-1}(x)\cap V_{i_0}|$ is maximum. Using the Pigeonhole
Principle we obtain  $b\geq (n-t)/(k-1)$ which gives $t\geq
n-(k-1)b$. Furthermore $b<n/k$ as there is only one vertex of $V(C)$
satisfying the negation ($f(i_0)\neq x$ already has  this
property). Now we define a mapping $f'\colon V(D)\to V(C)$ as  
  
      \[ f'(i)=\begin{cases}
      f(i) & \text{for}\;\; i\neq i_0\\
      x & \text{for}\;\; i=i_0.
      
      \end{cases}
      \]
  
 Since $f$ and $f'$ differ only at $i_0$ and $C$ is
 $D$-pointed, the function  $f'$ cannot be an acyclic
 homomorphism. We distinguish two cases.

\medskip
\textbf{Case I}: $x\not \in$ Range$(f)$.

In this case, the only reason that $f'$ is not an acyclic
homomorphism is that there must be a vertex $v\neq i_0$ in $V(D)$ such
that either  $f(v)f(i_0)\in A(C)$ but $f(v)x\not \in A(C)$
(and $vi_0\in A(D)$), 
or $f(i_0)f(v)\in A(C)$ but
$xf(v)\not \in A(C)$ 
(and $i_0v\in A(D)$). 
Without loss of
generality, assume that $f(v)f(i_0)\in A(C)$ but $f(v)x\not\in A(C)$ 
(and $vi_0\in A(D)$) occurs. 
We have 
$|\varphi^{-1}(f(v))\cap  V_v|\geq t\geq n-(k-1)b$, so we can  choose a set
$U\subseteq \varphi^{-1}(f(v))\cap V_v$ with $|U|=n-(k-1)b$. Then  there
must be at most $\min\{b,n^{\epsilon \ell}\}$ arcs from $U$ to
$\aA:=\varphi^{-1}(x)\cap V_{i_0}$ in $D'$; otherwise after passing from $D'$
to $D^*$, we have some arc(s) left between these two sets in $D^*$ and
since $\varphi$ is an acyclic homomorphism, $f(v)x\in A(C)$ which is a
contradiction. But the property just described is the rare property
articulated in Lemma~\ref{clm3}, and $D'$ was chosen not to enjoy it,
so Case~I leads to this contradiction. 
 
\medskip
\textbf{Case II}: $x\in$ Range$(f)$.

In this case, there are two potential reasons for $f'$ not to be an
acyclic homomorphism. The reason we explained in Case~I is still a
potential reason in the present case, and it similarly leads to a
contradiction. The other reason here is when $f'^{-1}(x)$ does not
induce an acyclic subdigraph in $D$. We proceed to show that
this also leads to a contradiction. 

We know that $\varphi^{-1}(x)\cap V_{i_0}\neq \varnothing$. Since
$x\in\text{Range}(f)$, we have $x=f(j_0)$ for some $j_0\in V(D)$ and
$j_0\neq i_0$. The reason for $j_0\neq i_0$ is that 
$f(i_0)\neq x=f(j_0)$. We show that in this case 
$i_0j_0, j_0i_0 \in A(D)$. Suppose to the contrary that this is
wrong. Without loss of 
generality, we assume that $i_0j_0\not\in A(D)$. First we claim that
there exists a vertex $p_0\in V(D)$ such that it has a different
situation with respect to $i_0$ and $j_0$ in the sense of adjacency
(like, for example, $p_0i_0\in A(D)$, but $p_0j_0\not\in A(D)$). For if
every $p_0\in V(D)$ that is adjacent to $i_0$ is also adjacent to
$j_0$ (preserving the directions), then we can define the mapping
$g\colon V(D)\to V(C)$ by $g(i)=f(i)$ for $i\neq i_0$ and
$g(i_0)=f(j_0)$. Then $f\neq g$ (but they differ only at $i_0$), and $g$ is
clearly an acyclic homomorphism; this contradicts the $D$-pointedness 
of $C$. We also claim that there exist $p_0\in V(D)$ and 
$v\in V(C)$ such that $f(p_0)=v$, the arc $p_0i_0\in A(D)$,  
$p_0j_0\not\in A(D)$, the arc $vf(i_0)\in A(C)$, and $vf(j_0)\not\in A(C)$. 
For if every $p_0\in V(D)$ with $p_0i_0\in A(D)$, $p_0j_0\not\in A(D)$ 
satisfies both $vf(i_0)\in A(C)$ and $vf(j_0)\in A(C)$, then we can
again define the mapping $g\colon V(D)\to V(C)$ by $g(i)=f(i)$ for 
$i\neq i_0$ and $g(i_0)=f(j_0)$, which again contradicts the fact
that $C$ is $D$-pointed.  

Thus let $p_0$, $v$ as above satisfy $f(p_0)=v$, the arc 
$p_0i_0\in A(D)$, $p_0j_0\not\in A(D)$,  the arc $vf(i_0)\in A(C)$, and  
$vf(j_0)\not\in A(C)$. The sets
$W_{p_0}:=\varphi^{-1}(v)\cap V_{p_0}$ and 
$\bB' :=\varphi^{-1}(f(j_0))\cap V_{i_0}$ satisfy 
$|W_{p_0}|\geq t=|\varphi^{-1}(f(i_0))\cap V_{i_0}|$
and $n/k\geq |\bB'|=b\geq (n-t)/(k-1)$. Hence, there exists a set
$\aA'\subseteq W_{p_0}$ such that $|\aA'|=n-|\bB'|(k-1)$ with the property that
there is no arc from $\aA'$ to $\bB'$ in $D$ (as $\varphi(\bB')=f(j_0)$ and
$\varphi(\aA')=v$ and $vf(j_0)\not\in A(C)$). However, this contradicts
Lemma~\ref{clm3}. Thus,  $i_0j_0, j_0i_0\in A(D)$. Using this important
fact, we proceed to show that (the second reason in) Case~II also leads to a
contradiction. 

The  definition of $f$ gives us  
$|\varphi^{-1}(f(j_0))\cap V_{j_0}|\geq n/k$. Since 
$n^{\epsilon \ell}\leq n^{1/4}\ll n/k$  we
can choose $\aA\subseteq \varphi^{-1}(f(j_0))\cap V_{j_0}$ with
$|\aA|=\lfloor n/2k \rfloor$ such that no arc of $M$ (the matching
defined at the start of Section~\ref{Sec:main_proof}) has an end
vertex in $\aA$.
Let $\bB=\{z\}\subset\varphi^{-1}(f(j_0))\cap V_{i_0}$. Since all arcs
of $M$ are independent, at most one arc of 
$M$ is incident with $z$. Since $\varphi(\aA\cup\bB)=\{x\}$ and $\varphi$
is an acyclic homomorphism, the subdigraph of $D^*$ induced by 
$\aA\cup\bB$ is acyclic. To show that this is unlikely, we first estimate the
expected number $N$ of ways to select a vertex $z\in V_{i_0}$ and a
subset $U\subseteq V_{j_0}$ of cardinality $\lfloor n/2k \rfloor$ so
that the subdigraph $H_{z,U}$ of $D^*$ they induce is acyclic and no
arc of $M$ is incident with a vertex in $U$. If $P_{z,U}$ denotes the
probability that $H_{z,U}$ is acyclic, then  
\begin{align}
    N\leq n\binom{n}{\lfloor n/2k\rfloor}P_{z,U}< n^nP_{z,U}.
    \label{firstupperboundforN}
\end{align} 

In order to bound $P_{z,U}$, we employ Chernoff's Inequality
(Theorem~\ref{Chernoffresult}). 
Let $\Omega$ be the set of all potential arcs in the subdigraph
$D'_{z,U}$, of $D_0$  induced by $\{z\}\cup U$. Each arc in $\Omega$
appears in $H_{z,U}$ with probability $p$. Let
$\tau>(2+\epsilon)/\epsilon$ be a fixed integer. We index (by positive
integers) those cycles of $D'_{z,U}$ that are of length $\tau+1$. For
$i\geq 1$, let $S_i$ be the arc set of the $i$th such cycle and $B_i$
be the event that the arcs in $S_i$ all appear (i.e., the cycle
determined by $S_i$ is present in $H_{z,U}$). Let the random variable
$Y$ count the $B_i$'s that occur. Since $P(Y=0)$ is an upper bound for
$P_{z,U}$, we can bound $P_{z,U}$ by bounding $P(Y=0)$. 
Using Theorem~\ref{Chernoffresult} with $\gamma=1$, we have 
\begin{align}
    P(Y=0)\leq P\big(|Y-E(Y)|\geq E(Y)\big)\leq 2e^{-E(Y)/3}.
    \label{boundonp(y=0)}
\end{align}

Since the arcs of $D'_{z,U}$ within $U$ are acyclically oriented, each
choice of $\tau$ vertices within $U$ determines exactly one potential
$(\tau+1)$-cycle. It follows that 
\begin{align}
    E(Y)=\binom{\lfloor n/2k \rfloor}{\tau}p^{\tau+1}>\Big(\frac{\lfloor n/2k \rfloor}{\tau}\Big)^{\tau}p^{\tau+1}>\frac{n^{\epsilon \tau+\epsilon-1}}{(4k\tau)^{\tau}}.
    \label{boundonmu}
\end{align}

\noindent 
Using (\ref{boundonp(y=0)}) and (\ref{boundonmu}) we find that
\begin{align*}
    P(Y=0)\leq 2e^{-n^{\epsilon \tau+\epsilon-1}/3(4k\tau)^{\tau}},
\end{align*}

\noindent 
and recalling our choice of $\tau$ (as exceeding 
$(2+\epsilon)/\epsilon$), we see that 
\begin{align}
    P(Y=0)\leq 2e^{-n^{1+2\epsilon }/3(4k\tau)^{\tau}}<e^{-n^{1+\epsilon}}.
    \label{finalupperboundforp(y=0)}
\end{align}

\noindent
Returning to (\ref{firstupperboundforN}), we have 
\begin{align*}
    N<n^nP_{z,U}<n^ne^{-n^{1+\epsilon}}=(ne^{-n^{\epsilon}})^n< e^{-n^{1+\epsilon}/2}.
\end{align*}

By Markov's Inequality, the probability that there exists such a set 
$\{z\}\cup U$ that induces an acyclic subdigraph  is less than
$e^{-n^{1+\epsilon}/2}$, which means it is unlikely as desired. 
  
Our discussion in Cases~I and II implies that $\varphi$ and $f \circ \psi$
have the same range.  It is now evident that $\varphi=f \circ \psi$, for
otherwise the same situation as in the proof that
$\text{Range}(\varphi)=\text{Range}(f\circ \psi)$ occurs and similarly leads to a
contradiction. Hence $\varphi=f \circ \psi$ as desired and therefore the
proof of part (iii) of Theorem~\ref{main} is complete.\hfill $\Box$
  
\subsection*{Acknowledgements}
The authors thank Liam Rafferty (who suggested this line of research
soon after completing \cite{Rafferty-thesis}) and Michael Morris (who
pointed out a gap in Case~II within the manuscript). Thanks also to
two referees ($A$ and $B$) for their careful reading and constructive comments.

%% References with bibTeX database:

%% \bibliographystyle{abbrv}
%% \bibliography{KayllParsa-2020}

%% References without bibTeX database:

% \begin{thebibliography}{00}

%% \bibitem must have the following form:
%%   \bibitem{key}...
%%

% \bibitem{}

% \end{thebibliography}

\end{document}